\documentclass[12pt,reqno]{amsart}
\usepackage{bm}
\theoremstyle{plain}
\newtheorem{thm}{Theorem}[section]
\newtheorem{cor}[thm]{Corollary}

\newtheorem{prop}[thm]{Proposition}
\theoremstyle{remark}
\newtheorem{rem}[thm]{Remark}
\theoremstyle{definition}
\newtheorem{defn}[thm]{Definition}

\newcommand{\rt}{{\bm{R}}}
\newcommand{\zt}{{\wit{Z}}}
\newcommand{\gnt}{{\bm{\gn}}}
\newcommand{\dq}{\partial_t}
\newcommand{\xt}{{\wit{X}}}
\newcommand{\yt}{{\wit{Y}}}
\newcommand{\dd}{\mathrm{d}}

\newcommand{\rl}{\mathbb{R}}

\newcommand{\wit}[1]{\widetilde{#1}}

\newcommand{\go}{\omega}

\newcommand{\gS}{\Sigma}
\newcommand{\gp}{\varphi}
\newcommand{\gn}{\nabla}
\newcommand{\eps}{\varepsilon}

\newcommand{\gG}{\Gamma}

\newcommand{\demi}{\frac{1}{2}}

\begin{document}
\title[The energy-momentum tensor as a second fundamental form]{The energy-momentum tensor as a second fundamental form}
\author{Bertrand Morel} 
\address{Bertrand Morel\newline
\indent Institut {\'E}lie Cartan\newline
\indent Universit\'e Henri Poincar\'e, Nancy I\newline
\indent B.P. 239\newline
\indent 54506 Vand\oe uvre-L{\`{e}}s-Nancy Cedex\newline
\indent France}
\email{morel@iecn.u-nancy.fr}



\begin{abstract}
We show that it is natural to consider the energy-mo\-men\-tum
tensor associated with a spinor field as the second fundamental form of an
isommetric immersion. In particular we give a generalization of the warped
product construction over a Riemannian manifold leading to this interpretation. Special
sections of the spinor bundle, generalizing the notion of Killing spinor, are
studied. First applications of such a construction are then given.
\end{abstract}

\maketitle

\section{Introduction}
In this paper, we show that there exists an analogy between the role
of the energy-momentum tensor associated with an eigenspinor field of the
Dirac operator and the role of the second fundamental form of a
hypersurface.

In \cite{Hi3}, O.~Hijazi proved that on a compact Riemannian spin manifold,
any eigenvalue $\lambda$ of the Dirac operator $D$, to which is attached an eigenspinor $\psi$, satisfies
\begin{equation}
\lambda^2\geq \inf_M(\frac{1}{4}S+|T^\psi|^2),\tag*{($\star$)}
\end{equation}
where $S$ is the scalar curvature of the manifold, and $T^\psi$ is the
energy-momentum tensor associated with $\psi$. This lower bound gives a
non-trivial information on the spectrum of $D$ without requiring $S$ to be
positive (compare with Friedrich's inequality \cite{Fri80}). Even though
the r.h.s. of ($\star$) depends on the eigenspinor $\psi$, we shall show
that, in the case of a hypersurface, the limiting case can be
geometrically interpreted.

Therefore, we begin by recalling basic facts regarding spin geometry of hypersurfaces, such
as the identification of the restriction of the spin bundle
of an ambient space with the spin bundle of a hypersurface, and the
spinorial Gau\ss{} formula.

We then prove that the extrinsic lower bound for the eigenvalue of the Dirac
operator on a compact hypersurface bounding a compact domain given in
\cite{HZM2} is related to the intrinsic estimate ($\star$). Equality
cases are characterized by the existence of special sections of the spinor
bundle, called \emph{$T$-Killing spinors}, which generalize the
notion of Killing spinors (these particular sections have been studied by E.C.~Kim and Th.~Friedrich in \cite{FK}).

Recall that complete simply connected Riemannian spin manifolds carrying
real Killing spinors are characterized in \cite{Bea}. For this, C.~B\"ar
proved that the usual cone constructed over such a manifold has to be Ricci flat. The geometry of
the cone being studied by S.~Gallot in \cite{Ga}, it suffices then to apply
the already known characterization of irreducible, simply connected
Riemannian spin manifolds carrying parallel spinors \cite{H},\cite{Wa}.

In the last three sections, we show that the above argument can be
generalized to Riemannian spin manifolds carrying $T$-Killing spinors
in the following way: 

We start by constructing a generalized warped product over a
Riemannian manifold $(M,g)$, deforming the initial metric in the
direction of the energy-momentum tensor $T^\psi$ associated with a $T$-Killing
spinor field $\psi$. Therefore, we can consider $(M,g)$ as a hypersurface of the
constructed manifold $\mathcal{M}$ whose second fondamental form is given
by $T^\psi$. In the case where this $2$-symmetric tensor is parallel, we then show how this
construction extends the one of the usual warped product. We finally prove that in the case of a manifold carrying a
$T$-Killing spinor field whose energy-momentum tensor is a projector, the
manifold $\mathcal{M}$ has to be Ricci flat (Theorem~\ref{thmprin}).  

\section{The restriction principle and $T$-Killing spinors}

\subsection{Restricting spinors to a hypersurface} 

Let $(N^{n+1},g)$ be an oriented Riemannian spin manifold of dimension $n+1$, with a fixed
spin structure. Denote by $\gS N$ the spinor bundle associated with this spin
structure. If $M^n$ is an oriented hypersurface isometrically immersed into
$N$, denote by $\nu$ its unit, globally defined, normal vector field. Then $M$ is endowed with a spin structure, canonically induced by
the one on $N$. Denote by $\gS M$ the corresponding spinor bundle.
Recall that the spinor bundle $\gS N$ splits into $$\gS
N=\gS^+N\oplus\gS^-N$$ where $\gS^\pm N$ is the $\pm 1$-eigenspace for the
action of the complex volume forme $\go_{n+1}=i^{[\frac{n+1}{2}]}\go$.
The following proposition is essential for what follows (see for example
\cite{Bar1}, \cite{BHMM}, \cite{Mor}, \cite{Trau1}):

\begin{prop}\label{3ident}
There exists an identification of $\gS N_{|M}$ (resp. $\gS^+ N_{|M}$) if
$n$ is even (resp. odd) with $\gS M$, which after
restriction to $M$, sends
every spinor field $\psi \in \gG(\gS N)$ to the spinor field denoted by
$\psi^* \in \gG(\gS M)$. Moreover, if $\underset{N}{\cdot}$ (resp.
$\cdot$) stands for Clifford multiplication on $\gS N$ (resp. $\gS M$),
then one has
\begin{equation}
  \label{3eq:clifmul}
  (X\underset{N}{\cdot}\nu\underset{N}{\cdot}\psi)^*=X\cdot\psi^*\;,
\end{equation}
for any vector field $X$ tangent to $M$.
\end{prop}

Another important formula is the well-known spinorial Gauss formula: if
$\gn^N$ and $\gn$ stand for the covariant derivatives on $\gG(\gS N)$ and
$\gG(\gS M)$  respectively, then, for all $X \in TM$ and $\psi \in
\gG(\gS N)$
\begin{equation}
  \label{3eq:gauss}
  (\gn^N_X\psi)^*=\gn_X\psi^*+\demi h(X)\cdot\psi^*,
\end{equation}
where $h$ is the second fundamental form of the immersion $M\hookrightarrow
N$ viewed as a symmetric endomorphism of the tangent bundle of $M$.

Recall that the ambient spinor bundle $\gS N$ can be endowed  with a Hermitian inner product $(\;,)_N$ for which Clifford
multiplication by any vector tangent to $N$ is skew-symmetric. This product
induces another Hermitian inner product on $\gS M$, denoted by $(\;,)$ making the
identification of Proposition \ref{3ident} an isometry. Now, relation
(\ref{3eq:clifmul}) shows that Clifford
multiplication by any vector tangent to $M$ is skew-symmetric with respect
to $(\;,)$. 

\subsection{On Hijazi's inequality involving the energy-momentum tensor}

We now discuss the role of the energy-momentum tensor associated with a
special section of the spinor bundle when it is involved in lower bounds
for the first eigenvalue of the Dirac operator.
Recall the following estimate

\begin{thm}[Hijazi, \cite{Hi3}]\label{3hijaz} On a compact Riemannian spin manifold
  $(M^n,g)$, any eigenvalue $\lambda$ of the Dirac operator to which is
  attached an eigenspinor $\psi$, satisfies
  \begin{equation}
    \label{3eq:hija}
    \lambda^2\geq \inf_M(\frac{1}{4}S+|T^\psi|^2),
  \end{equation}
where $S$ is the scalar curvature of $(M^n,g)$, and $T^\psi$ the field of
symmetric endomorphisms defined on the complement of the set of zeros of
$\psi$ by
$$T^\psi(X,Y)=\demi\Re(X\cdot\gn_{Y}\psi+Y\cdot\gn_{X}\psi,\frac{\psi}{|\psi|^2})\;.$$

\end{thm}

By the Cauchy-Schwarz inequality, one can easily check that this lower bound improves the one obtained by Th.~Friedrich \cite{Fri80}
$$\frac{n}{4(n-1)}\inf_M S\;,$$ and in particular, it may give a non-trivial
lower bound in the case where $\inf_M S$ is non-positive whereas
Friedrich's inequality requires $S$ to be positive.

If equality holds in (\ref{3eq:hija}), then $\psi$ has constant length,
\begin{equation}
  \label{3eq:trac1}
  \mathrm{tr}(T^\psi)^2=\frac{1}{4}S+|T^\psi|^2=\text{constant} 
\end{equation}and
\begin{equation}
  \label{3eq:trac2}
  \mathrm{div}(T^\psi)=0\;.
\end{equation}

Interpreting the energy-momentum tensor associated with a
spinor field as a second fundamental form (see also \cite{FK}) leads to the following proposition and remarks:  

\begin{prop}\label{3prorelat} Let $M^n\hookrightarrow (N^{n+1},g)$ be any compact oriented hypersurface
isometrically immersed in an oriented Riemannian spin manifold
$(N^{n+1},g)$, with constant mean curvature $H$ and second fundamental form
$h$. Assume $(N^{n+1},g)$
admits a parallel spinor field, then $M$
satisfies the equality case in (\ref{3eq:hija}).

Moreover, the energy-momentum tensor associated with the spinor field $\psi$
in Theorem \ref{3hijaz} satisfies $$2T^\psi=h\;.$$
\end{prop}

\begin{proof}
Let $\Phi$ be any parallel spinor field on $N$. Then Proposition \ref{3ident} and Gau{\ss} formula (\ref{3eq:gauss}) yield 
\begin{equation}\label{3emspin}
\gn_X\psi+\demi h(X)\cdot\psi=0\;,
\end{equation}
 with $\psi:=\Phi^*$.  Let  $(e_1,\ldots, e_n)$ be a positively oriented local orthonormal basis of
$\gG(TM)$. Then for $j=1,\ldots,
 n$ we have $$\gn_{e_j}\psi=-\sum_{k=1}^n
\demi h_{jk}e_k\cdot\psi\;.$$
Taking Clifford multiplication by $e_i$ and the scalar product with $\psi$,
we get
$$\Re(e_i\cdot\gn_{e_j}\psi,\psi)=-\sum_{k=1}^2 \demi h_{jk}\Re(e_i\cdot
e_k\cdot\psi,\psi)\;.$$
Since $\Re(e_i\cdot
e_k\cdot\psi,\psi)=-\delta_{ik}|\psi|^2$, it follows, by the symmetry of $h$
$$\Re(e_i\cdot\gn_{e_j}\gp+e_j\cdot\gn_{e_i}\psi,\psi)=h_{ij}|\psi|^2\;.$$
Therefore, $2T^\psi=h$. Moreover, $N$ has to be Ricci-flat,
and the Gau{\ss} Equation yields  $$\mathrm{tr}(2T^\psi)^2=n^2
H^2=S+|2T^\psi|^2=\text{constant}\;.$$ Tracing Equation
(\ref{3emspin}) completes the proof. 
\end{proof} 

\begin{rem}\label{rem1}
Since the Gau{\ss} equation implies $$\frac{n^2}{4} H^2=
\frac{1}{4}S+|T^\psi|^2=\text{constant}\;,$$ we point out that we are exactly
in the equality case of the extrinsic estimate given in \cite{HZM2}
(Theorem 6), as expected. 
\end{rem}

\begin{rem}
Conversely, we can not conclude that if $(M,g)$, $\dim M\geq 3$, satisfies the equality
case in (\ref{3eq:hija}), there exists an isometric immersion of $M$ as a hypersurface with constant mean curvature 
into a Riemmanian spin manifold with parallel
spinors. However, this is true if $\dim M=2$ (see \cite{Fr}).  
\end{rem}

\subsection{On $T$-Killing spinors}

We now give basic properties of $T$-Killing spinors,
which satisfy equality case in (\ref{3eq:hija}). They are generalizations of Killing spinors, and have been
studied by Th.~Friedrich and E.C.~Kim in \cite{FK}.

Let $(M^n,g)$ be a $n$-dimensional Riemannian spin manifold. A $T$-Killing
spinor field $\psi \in \Gamma(\Sigma M)$ is a spinor field which satisfies
\begin{equation}
  \label{eq:emspin}
\gn_X\psi=-T^\psi(X)\cdot\psi\;,\quad \forall X\in TM  
\end{equation}
and
$$\mathrm{tr}(T^\psi)=\; \text{constant}, $$ where $T^\psi$ is the
energy-momentum tensor associated with $\psi$, defined as in Theorem
\ref{3hijaz}. It is easy to see that such a spinor field is an eigenspinor
for the Dirac operator and since it satisfies the equality case in
(\ref{3eq:hija}), it also satisfies Equations (\ref{3eq:trac1}) and
(\ref{3eq:trac2}). Moreover $\psi$ has constant length. 

Indeed, computing the action of the spinorial curvature tensor $\mathcal{R}$ on the
spinor $\psi$, we see that necessarily, for all vector fields $X, Y\in \Gamma(TM)$,
\begin{eqnarray}
  \label{3eq:TcondR}
  \mathcal{R}(X,Y)\psi&=&\gn_X\gn_Y\psi-\gn_Y\gn_X\psi-\gn_{[X,Y]}\psi\nonumber\\
  &=&\Big ( T^\psi(Y)\cdot T^\psi(X)-T^\psi(X)\cdot T^\psi(Y)\Big
  )\cdot\psi\\
  &&+\Big ( (\gn_Y T^\psi)(X)-(\gn_X T^\psi)(Y)\Big
  )\cdot\psi\nonumber
\end{eqnarray}
and
\begin{eqnarray}
  \label{3eq:RicTpsi}
  \mathrm{Ric}(X)\cdot\psi&=&2\sum_{i=1}^n e_i\cdot\mathcal{R}(e_i,X)\psi\nonumber\\
&=&4\mathrm{tr}(T^\psi)T^\psi(X)\cdot\psi-4T^2_\psi(X)\cdot\psi-2\sum_{i=1}^{n}e_i\cdot(\gn_{e_i}T^\psi)(X)\cdot\psi\;.
\end{eqnarray}

Recall the following definition

\begin{defn} A symmetric $2$-tensor $T\in S^2(M)$ is called a \emph{Codazzi
    tensor} if it satisfies the Codazzi-Mainardi equation, i.e. 
  \begin{equation}
    \label{eq:cod}
 (\gn_XT)(Y)=(\gn_YT)(X)\qquad\forall X,Y\in\gG(TM)\;,   
\end{equation}
($T$ being viewed in this formula via the metric $g$ as a symmetric endomorphism of the
tangent bundle).
\end{defn}

To every nowhere vanishing spinor field $\psi$, we can associate the real vector field $V_\psi$
defined by
$$g(V_\psi,X)=i(\psi,T^\psi(X)\cdot\psi)\;,\qquad \forall X \in \Gamma(TM)\;.$$

Note that $T$-Killing spinors have the following property:

\begin{prop}\label{3pprroo}
Let $\psi$ be a $T$-Killing spinor field. If $T^\psi$ is a Codazzi tensor,
then the associated vector field $V_\psi$ is a Killing vector field.
\end{prop}
\begin{proof} For all $X, Y \in TM$, we compute
  \begin{eqnarray*}
  g(\gn_XV_\psi,Y)&\!=&\!Xg(V_\psi,Y)- g(V_\psi,\gn_XY)\\
&\!=&\!i(\gn_X\psi,T^\psi(Y)\cdot\psi)+i(\psi,\gn_X(T^\psi(Y))\cdot\psi)\\&&+i(\psi,T^\psi(Y)\cdot\gn_X\psi)-i(\psi,T^\psi(\gn_XY)\cdot\psi)\end{eqnarray*}
and since Clifford multiplication by vector fields is skew-symmetric with
respect to $(.\,,.)$, we have
$$ g(\gn_XV_\psi,Y)=i(T^\psi(Y)\cdot\! T^\psi(X)\!\cdot\psi,\psi)-i(T^\psi(X)\cdot
 T^\psi(Y)\cdot\psi,\psi)+i(\psi,\gn_X(T^\psi)(Y)\cdot\psi) $$
which is clearly skew-symmetric if \eqref{eq:cod} holds for $T^\psi$.
\end{proof}

\begin{rem}\label{rem2}
A spinor field satisfying Equation \eqref{eq:emspin} and whose associated
energy-mo\-men\-tum tensor $T^\psi$ is a Codazzi tensor is called a
\emph{Codazzi Energy-Momentum spinor}. This notion generalizes the notion
of Killing spinors (see \cite{MoCEM} for a study of these particular spinor
fields). It is known that given a Riemannian manifold $(M,g)$ of dimension
$3$, the existence of a non trivial Codazzi Energy-Momentum spinor field
is equivalent to the existence of an isometric immersion of the universal
covering of $M$ into the Euclidean $4$-dimensional space (see \cite{MoS3}).
\end{rem}

\section{Generalized warped product over a Riemannian manifold}

Proposition \ref{3prorelat} and Remarks \ref{rem1} and \ref{rem2} show that
it is natural to consider the energy-momentum tensor associated with a spinor
field satisfying \eqref{eq:emspin} as the second fondamental form of an
apropriate immersion. In this section, we give a natural construction of a
generalized warped product over a Riemannian manifold.

Recall that if $\pi: M\rightarrow N$ is a smooth map between two
manifolds, then two vector fields $X\in\Gamma(TM)$ and $Y\in\Gamma(TN)$ are said to be
$\pi$-related ($X\overset{\pi}\sim Y$) if $\dd
\pi(X_p)=Y_{\pi(p)}$ at all point $p$ in $M$. Two vector fields
$X\in\Gamma(TM)$ and $Y\in\Gamma(TN)$ are $\pi$-related if and only if
$X(f\circ\pi)=Y(f)\circ\pi$ for all smooth function $f$ defined on
$N$. Moreover, it is a classical fact that  $X_1\overset{\pi}\sim X_2$ and
$Y_1\overset{\pi}\sim Y_2$ imply $[X_1,Y_1]\overset{\pi}\sim [X_2,Y_2]$.

Now let $(M^n,g)$ be a compact Riemannian spin manifold. We denote by $\gn$ the
Levi-Civita connection on $M$.
Let $h\in S^2(M)$ be a symmetric $2$-tensor, and $f:I\rightarrow \rl$ a
smooth function, with $f(1)=0$ and $I=]1-\eps,1+\eps[\subset \rl$, such that the metric 
\begin{equation}\label{3metric}
g_t=g+ f(t)h
\end{equation} is well-defined on $M$ for all $t\in I$.
We endow  $M\times I$ with the metric
\begin{equation}
  \label{3metrich}
  \mathbf{g}=\pi_M^*(g_t)\oplus\pi_I^*(\dd^2 t)\, ,
\end{equation}
where $\pi_M$ and $\pi_I$ are respectively the first and second canonical
projections of $M\times I$ on $M$ and $I$.
We denote by $\mathcal{M}=(M\times I,\mathbf{g})$ the generalized Riemannian
warped product obtained by this construction.

Note that for $f(t)=t^2-1$ and $h=g$, this construction corresponds to the
usual cone over $M$, the metric $\mathbf{g}$ being then defined by
$\mathbf{g}=t^2g+\dd^2 t$ (see \cite{Ga}). In the following, we will call
$\mathcal{M}$ a {\it generalized cone} if by definition $f(t)=t^2-1$

The lift of a vector field $X\in\Gamma(TM)$ to
$\mathcal{M}$ will be called a {\it horizontal} vector field. By
definition, it is the only vector field, $\xt\in\Gamma(T\mathcal{M})$, such
that $\dd\pi_M(\xt)=X$ and  $\dd\pi_I(\xt)=0$. The set of horizontal vector
fields will be denoted by $\mathcal{H}(\mathcal{M})$. In the following, $\dq$ will
stand for the unit vector field spanning $\Gamma(TI)$, as well as its
lift to $\Gamma(T\mathcal{M})$. Then we have the following classical
proposition (see \cite{Onil} for example):

\begin{prop}
If $\xt,\yt \in \mathcal{H}(\mathcal{M})$ are horizontal vector fields, then
$$[\xt,\yt]=\wit{[X,Y]},\quad \mbox{and}\quad [\xt,\dq]=0.$$
If $q\,: I\rightarrow \rl$ is a smooth function, then the gradient of
$\wit{q}=q\circ\pi_I\,:\mathcal{M}\rightarrow \rl$ is the lift of the
gradient of $q$.
\end{prop} 

For all  $X\in\Gamma(TM)$ and $t\in I$ we denote by $H_t(X)$ the vector
field defined by  $$g_t(H_t(X),Y)=h(X,Y),\quad\forall Y\in\Gamma(TM).$$

\begin{prop}\label{3prop1}
If $\,\gnt$ stands for the Levi-Civita connection of the generalized
Riemannian warped product $\mathcal{M}$, and $\gn^t$ for the Levi-Civita connection
of $(M,g_t)$, then for all $X,Y \in \Gamma(TM)$, we have
the following generalized O'Neill formulas:
\begin{eqnarray}
\gnt_{\dq}\dq&=&0\;,\label{nabla1}\\
\gnt_{\xt}\dq=\gnt_{\dq}\xt&=&\frac{f'(t)}{2}\wit{H_t(X)}\;,\label{nabla2}\\
\gnt_{\xt}\yt&=&\wit{\gn^t_XY}-\frac{f'(t)}{2}h(X,Y)\dq\;.\label{nabla3}
\end{eqnarray}
\end{prop}

\begin{proof} Equality \eqref{nabla1} is trivial. Hence, note that since
  $\mathbf{g}(\xt,\dq)=0$, we have $$\mathbf{g}(\gnt_{\dq}\xt,\dq)=0\;.$$ On the other hand, $\mathbf{g}(\dq,\dq)=1$
  implies $\mathbf{g}(\gnt_{\xt}\dq,\dq)=0$. By the Koszul formula, we get
\begin{eqnarray*}
2\mathbf{g}(\gnt_{\dq}\xt,\yt)&=&2\mathbf{g}(\gnt_{\xt}\dq,\yt)=\dq(\mathbf{g}(\xt,\yt))=\dq(g_t(X,Y))\\
&=&\dq(g(X,Y))+\dq(f(t)h(X,Y))=f'(t)h(X,Y)\\
&=&f'(t)g_t(H_t(X),Y)=f'(t)\mathbf{g}(\wit{H_t(X)},\yt),
\end{eqnarray*}which proves \eqref{nabla2}. Again, by the Koszul formula
\begin{eqnarray*}
2\mathbf{g}(\gnt_{\xt}\yt,\dq)&=&-\dq(\mathbf{g}(\xt,\yt))=-\dq(g_t(X,Y))\\
&=&-\dq(g(X,Y))-\dq(f(t)h(X,Y))=-f'(t)h(X,Y)\;,\end{eqnarray*}and
\begin{eqnarray*}
2\mathbf{g}(\gnt_{\xt}\yt,\zt)&=&\xt g_t(Y,Z)+\yt g_t(Z,X)-\zt g_t(X,Y)\\&&-g_t(X,[Y,Z])+g_t(Y,[Z,X])+g_t(Z,[X,Y])\\&=&2g_t(\gn^t_XY,Z)=2\mathbf{g}(\wit{\gn^t_XY},\zt)\;.\end{eqnarray*}
Therefore the proof of \eqref{nabla3} is completed.\end{proof}
\begin{prop}\label{3prop2}
For all $X,Y,Z\in\Gamma(TM)$, we have
\begin{eqnarray}
\gnt_{\xt}\gnt_{\dq}\dq=\gnt_{\dq}\gnt_{\dq}\dq&=&0\label{riem1}\\
\gnt_{\dq}\gnt_{\xt}\dq=\gnt_{\dq}\gnt_{\dq}\xt&=&\frac{f''(t)}{2}\wit{H_t(X)}-\frac{(f'(t))^2}{4}\wit{H_t(H_t(X))}\label{riem2}\\
\gnt_{\xt}\gnt_{\dq}\yt=\gnt_{\xt}\gnt_{\yt}\dq&=&\frac{f'(t)}{2}\wit{\gn^t_XH_t(Y)}-\frac{(f'(t))^2}{4}h(X,H_t(Y))\dq\label{riem3}\\
\gnt_{\dq}\gnt_{\xt}\yt&=&\frac{f'(t)}{2}\wit{A_t(X,Y)}-\frac{f'(t)}{2}\wit{H_t(\gn^t_XY)}-\frac{f''(t)}{2}h(X,Y)\dq\label{riem4}\\
\gnt_{\xt}\gnt_{\yt}\zt&=&\wit{\gn^t_X\gn^t_YZ}-\frac{f'(t)}{2}h(X,\gn^t_YZ)\dq\nonumber\\&&-\frac{f'(t)}{2}Xh(Y,Z)\dq-\frac{(f'(t))^2}{4}h(Y,Z)\wit{H_t(X)}\;,\label{riem5}\end{eqnarray}
where $A_t(X,Y)$ is the vector field defined on $M$ by \begin{eqnarray*}g_t(A_t(X,Y),Z)&=& Xh(Y,Z)+ Yh(Z,X)-Zh(X,Y)\\&&-h(X,[Y,Z])+h(Y,[Z,X])+h(Z,[X,Y])\;.\end{eqnarray*}\end{prop}
\begin{proof}
Equality \eqref{riem1} is trivial. Since $\mathbf{g}(\gnt_{\xt}\dq,\dq)=0$, we have $\mathbf{g}(\gnt_{\dq}\gnt_{\xt}\dq,\dq)=0.$ Moreover,
\begin{eqnarray*}
 \mathbf{g}(\gnt_{\dq}\gnt_{\xt}\dq,\yt)&=& \dq \mathbf{g}(\gnt_{\xt}\dq,\yt)-
 \mathbf{g}(\gnt_{\xt}\dq,\gnt_{\dq}\yt)\\&=&\dq
 \Big(\frac{f'(t)}{2}\mathbf{g}(\wit{H_t(X)},\yt)\Big)-\frac{(f'(t))^2}{4}
 \mathbf{g}(\wit{H_t(X)},\wit{H_t(Y)})\\&=&\dq
 \Big(\frac{f'(t)}{2}g_t(H_t(X),Y)\Big)-\frac{(f'(t))^2}{4}
 g_t(H_t(X),H_t(Y))\\&=&\frac{f''(t)}{2}g_t(H_t(X),Y)-\frac{(f'(t))^2}{4}
 g_t(H_t(H_t(X)),Y)\\&=&\mathbf{g}(\frac{f''(t)}{2}\wit{H_t(X)}-\frac{(f'(t))^2}{4}\wit{H_t(H_t(X))},Y)\;,
\end{eqnarray*} which proves \eqref{riem2}. On the other hand,
\begin{eqnarray*}
  \mathbf{g}(\gnt_{\xt}\gnt_{\dq}\yt,\zt)&=&\xt \mathbf{g}(\gnt_{\dq}\yt,\zt)- \mathbf{g}(\gnt_{\dq}\yt,\gnt_{\xt}\zt)\\
&=&\frac{f'(t)}{2}\Big(\xt g_t(H_t(Y),Z)-g_t(H_t(Y),\gn^t_XZ)\Big )\\&=&\frac{f'(t)}{2}\mathbf{g}(\wit{\gn^t_XH_t(Y)},\zt)\;,
\end{eqnarray*}and
\begin{eqnarray*}
   \mathbf{g}(\gnt_{\xt}\gnt_{\dq}\yt,\dq)&=&\xt \mathbf{g}(\gnt_{\dq}\yt,\dq)- \mathbf{g}(\gnt_{\dq}\yt,\gnt_{\xt}\dq)\\
&=&-\frac{(f'(t))^2}{4}g_t(H_t(Y),H_t(X))\\&=&-\frac{(f'(t))^2}{4}h(H_t(X),Y)\;.
\end{eqnarray*} Therefore, we proved \eqref{riem3}. For \eqref{riem4}, we
show that 
\begin{eqnarray*}
 \mathbf{g}(\gnt_{\dq}\gnt_{\xt}\yt,\dq)&=&\dq \mathbf{g}(\gnt_{\xt}\yt,\dq)- \mathbf{g}(\gnt_{\xt}\yt,\gnt_{\dq}\dq)\\&=&-\dq(\frac{f'(t)}{2}h(X,Y))\\&=&-\frac{f''(t)}{2}h(X,Y)
\end{eqnarray*}and
\begin{eqnarray*}
 \mathbf{g}(\gnt_{\dq}\gnt_{\xt}\yt,\zt)&=& \dq \mathbf{g}(\gnt_{\xt}\yt,\zt)-\mathbf{g}(\gnt_{\xt}\yt,\gnt_{\dq}\zt)\\&=& \dq \mathbf{g}(\wit{\gn^t_XY},\zt)-\frac{f'(t)}{2}\mathbf{g}(\wit{\gn^t_XY},\wit{H_t(Z)})\\&=&\frac{f'(t)}{2}\mathbf{g}(\wit{A_t(X,Y)},\zt)-\frac{f'(t)}{2}\mathbf{g}(\wit{H_t(\gn^t_XY)},\zt)\;.
\end{eqnarray*} On the other hand, we have
\begin{eqnarray*}
  \mathbf{g}(\gnt_{\xt}\gnt_{\yt}\zt,\dq)&=& \xt \mathbf{g}(\gnt_{\yt}\zt,\dq)- \mathbf{g}(\gnt_{\yt}\zt,\gnt_{\xt}\dq)\\&=&-\frac{f'(t)}{2}\Big( Xh(Y,Z)+ g_t(\gn^t_YZ,H_t(X))\Big )\\&=&-\frac{f'(t)}{2}\Big( Xh(Y,Z)+ h(\gn^t_YZ,X)\Big )
\end{eqnarray*} and
\begin{eqnarray*}
  \mathbf{g}(\gnt_{\xt}\gnt_{\yt}\zt,\wit{V})&=& \xt \mathbf{g}(\gnt_{\yt}\zt,\wit{V})- \mathbf{g}(\gnt_{\yt}\zt,\gnt_{\xt}\wit{V})\\&=&\xt g_t(\gn^t_YZ,V)-g_t(\gn^t_YZ,\gn^t_XV)-\frac{(f'(t))^2}{4}h(Y,Z)h(X,V)\\&=&g_t(\gn^t_X\gn^t_YZ,V)-\frac{(f'(t))^2}{4}h(Y,Z)g_t(H_t(X),V)\\&=&\mathbf{g}(\wit{\gn^t_X\gn^t_YZ},\wit{V})-\frac{(f'(t))^2}{4}h(Y,Z)\mathbf{g}(\wit{H_t(X)},\wit{V})\;.
\end{eqnarray*}
Therefore the proof of \eqref{riem5} is completed.
\end{proof}
\begin{prop}\label{3prop3}
If $\rt$ and $R^t$ denote respectively the Riemann curvature tensor of
$\mathcal{M}$ and of $(M,g_t)$, then we have the following relations:
\begin{eqnarray*}
\rt(\dq,\dq)\dq&=&\rt(\dq,\dq)\xt=0\\
\rt(\xt,\dq)\dq&=&-\frac{f''(t)}{2}\wit{H_t(X)}+\frac{(f'(t))^2}{4}\wit{H_t(H_t(X))}\\
\rt(\xt,\dq)\yt&=&\frac{f'(t)}{2}\Big
(\wit{\gn^t_XH_t(Y)}+\wit{H_t(\gn^t_XY)}-\wit{A_t(X,Y)}\Big )\\
&&+\frac{f''(t)}{2}h(X,Y)\dq-\frac{(f'(t))^2}{4}h(X,H_t(Y))\dq\\
\rt(\xt,\yt)\dq&=&\frac{f'(t)}{2}\Big[\wit{\gn^t_X(H_t)(Y)}-\wit{\gn^t_Y(H_t)(X)}\Big ]\\
\rt(\xt,\yt)\zt&=&\wit{R^t(X,Y)Z}+\frac{f'(t)}{2}\Big ((\gn_Yh)(X,Z)-(\gn_Xh)(Y,Z)\Big )\dq\\&&+\frac{(f'(t))^2}{4}\Big (h(X,Z)\wit{H_t(Y)}-h(Y,Z)\wit{H_t(X)}\Big )\;.
\end{eqnarray*}\end{prop}
\begin{proof} Straightforward from the preceding proposition.\end{proof}

\begin{rem}\label{3remoneil}
Let $q$ be a positive smooth function on $I$. Define $f:=q^2-1$
and $h=g$. Then $g_t=q^2(t)g$ and $\mathcal{M}$ is the usual warped
product constructed over $M$. We have $f'=2qq'$ and $f''=2(qq''+(q')^2)$, and since
$$q^2(t)g(H_t(X),Y)=g_t(H_t(X),Y)=h(X,Y)=g(X,Y)\;,$$ we get
$$H_t(X)=\frac{1}{q^2(t)}X\;.$$ Therefore, it is straightforward to see that Propositions \ref{3prop1}, \ref{3prop2}
and \ref{3prop3} correspond precisely to O'Neill Formulas (\cite{Onil} p.\,206 and 210). 
\end{rem} 

In the following, take $f(t)=t^2-1$. If $X,Y,Z,W\in \gG(T(M\times_f I)_{|M})$ are vector fields tangent to $M$
and $\dq$ the unit normal vector field on $M$, then, Proposition \ref{3prop3} yields
\begin{eqnarray*}
\mathbf{g}(\rt(X,\dq)\dq,Y)&=&h(H(X),Y)-h(X,Y)\\\mathbf{g}(\rt(X,Y)\dq,Z)&=&(\gn_Xh)(Y,Z)-(\gn_Yh)(X,Z)\\\mathbf{g}(\rt(X,Y)Z,W)&=&g(R(X,Y)Z,W)+h(X,Z)h(Y,W)-h(Y,Z)h(X,W)\;.\end{eqnarray*}
As expected, the construction of $\mathcal{M}$ allows the interpretation of $h$ as the second
fundamental form of the hypersurface $M\times\{1\}\subset \mathcal{M}$.

\section{$T$-Killing spinors with parallel Energy-Momentum tensor}

Assume now that $(M^n,g)$ admits a non trivial $T$-Killing
spinor field $\psi$ with parallel energy-momentum tensor $T^\psi$. Consider
the generalized cone $\mathcal{M}$ over $M$, with $h=2T^\psi$. 

Back to the construction of the cone, given $X,Y\in TM$, it is straigthforward to see that, 
$$\gn^t_XY=\gn_XY+(t^2-1)B^t(X,Y)$$ where the vector valued $2$-symmetric
tensor $B^t$ is defined by
$$g_t(B^t(X,Y),Z)=(\gn_Xh)(Y,Z)+(\gn_Yh)(Z,X)-(\gn_Zh)(X,Y)\;,\quad\forall
Z\in TM\;.$$

Therefore, Propositions \ref{3prop1}, \ref{3prop3} and formula
(\ref{3eq:RicTpsi}) yield
\begin{prop}\label{3prop1par}
  Assume that $(M^n,g)$ admits a non trivial $T$-Killing
spinor field $\psi$ with parallel energy-momentum tensor $T^\psi$. Then
the Levi-Civita connection and the Riemann curvature tensor of the generalized cone $\mathcal{M}$ over $M$, with $h=2T^\psi$ satisfy
\begin{eqnarray*}
\gnt_{\dq}\dq&=&0\\
\gnt_{\xt}\dq=\gnt_{\dq}\xt&=&t\,\wit{H_t(X)}\\
\gnt_{\xt}\yt&=&\wit{\gn_XY}-t\,h(X,Y)\dq\;,\end{eqnarray*}
and
\begin{eqnarray*}
\rt(\dq,\dq)\dq&=&\rt(\dq,\dq)\xt=\rt(\xt,\yt)\dq=0\\
\rt(\xt,\dq)\dq&=&-\wit{H_t(X)}+t^2\wit{H_t(H_t(X))}\\
\rt(\xt,\dq)\yt&=&h(X,Y)\dq-t^2h(X,H_t(Y))\dq\\
\rt(\xt,\yt)\zt&=&\wit{R(X,Y)Z}+t^2\Big (h(X,Z)\wit{H_t(Y)}-h(Y,Z)\wit{H_t(X)}\Big )\;.
\end{eqnarray*}
Moreover, the Ricci tensor of $(M^n,g)$, viewed as a field of symmetric endomorphism, satisfies for all $X\in \Gamma(TM)$
$$ \mathrm{Ric}(X)=\mathrm{tr}(H)H(X)-H^2(X)$$ where $H$ stands for $H_1$
with the above notations.
\end{prop}

With the help of this proposition, we can give an explicit relation between the Ricci tensor of the
generalized cone $\mathcal{M}$ and the Ricci tensor of $M$. For this, we
have to define a canonical local oriented orthonormal basis of $T\mathcal{M}$. At a
fixed $t \in I$, define the symmetric positive endomorphism $G_t$ of $TM$
by
$$g_t(X,Y)=g(G_t(X),Y)\;,\quad \forall X,Y\in TM\;.$$
Now, given a local oriented orthonormal basis $(e_1,\ldots ,e_n)$ of $TM$,
we get the local oriented orthonormal basis $(E_1,\ldots, E_n,\dq)$ of
$T\mathcal{M}$ by defining, at each point $(x,t)\in\mathcal{M}$, the vector
$E_i(x,t)$
as the lift of the vector $$e^t_i(x):= G_t^{-\demi}(e_i(x)) \in T_xM$$ (where for
all positive symmetric endomorphism $B$, $B^{\demi}$ stands for its
positive square root).

Since $g_t=g+(t^2-1)h$, it is easy to see that $G_t=\mathrm{Id}+(t^2-1)H$
for all $t\in I$. Hence, at a fixed point $t \in I$, $G_t$ is parallel
since we assumed that $H$ is parallel.

We then have the following

\begin{cor}\label{3corRic} Denote the Ricci tensor of $\mathcal{M}$ by $\wit{Ric}$, then
  for all horizontal vector fields $\xt,\yt \in \mathcal{H}(\mathcal{M})$, we
  have at a fixed point $(x,t)\in \mathcal{M}$
  \begin{eqnarray*}
    \wit{Ric}(\xt,\dq)&=&0\\ 
\wit{Ric}(\dq,\dq) &=&t^2\sum_{i=1}^n h(e_i^t,H_t(e_1^t))-\sum_{i=1}^n h(e_i^t,e_i^t)\\
\wit{Ric}(\xt,\yt) &=&\mathrm{tr}(H)h(X,Y)-h(H(X),Y)-t^2\sum_{i=1}^n
h(e_i^t,e_i^t)h(X,Y)+t^2 h(H_t(X),Y)\;.
  \end{eqnarray*}
\end{cor}

\begin{proof} Straightforward from Proposition \ref{3prop1par}. Remark that
  since $G_t$ is parallel, we have $$g_t(R(e_i^t,X)e_i^t,Y)=g(R(e_i,X)e_i,Y)\;.$$
\end{proof}

\begin{rem}\label{remmer}
As in Remark \ref{3remoneil}, assuming $h=g$, then $g_t=t^2g$ and
$\mathcal{M}$ is the usual cone over $M$. We get
$$H_t(X)=\frac{1}{t^2}X\qquad\text{and}\qquad e_i^t=\frac{1}{t}e_i\;.$$
Since $2T^\psi=h$, the spinor field $\psi\in\Gamma(\Sigma M)$ is actually a
real Killing spinor field with Killing number $\demi$ and we recover C.~B\"ar's
result which says that $\mathcal{M}$ has to be Ricci flat (\cite{Bea}).
\end{rem}

\section{The case of a projector}

If $M$ carries a $k$-dimensional
parallel smooth distribution $\mathcal{K}$, we can define on $M$ a parallel
field of symmetric endomorphisms $H$ satisfying $H^2=H$ as the projector on
$\mathcal{K}^\perp$, the orthogonal of $\mathcal{K}$. We can also consider the generalized cone
$\mathcal{M}$ over $M$ constructed with $H$. Then, if
$(e_1,\ldots , e_n)$ is an oriented orthonormal basis of $TM$ 
 such that the vectors $e_1,\ldots ,e_k$ span $\mathcal{K}$, we get 
$$H_t(X)=\begin{cases}0&\text{if }X\in\mathcal{K}\\\frac{1}{t^2}X&\text{if }X\in\mathcal{K}^{\perp}\end{cases}$$
and
$$e_i^t=\begin{cases}e_i&\text{if }1\leq i\leq k\\
\frac{1}{t}e_i&\text{if }k+1\leq i\leq n\;.\end{cases}$$

Hence, Corollary \ref{3corRic} implies the following

\begin{thm}\label{thmprin}
Let $(M^n,g)$ be a compact Riemannian spin manifold admitting
a smooth parallel distribution $\mathcal{K}$ and let $\psi$ be a
non-trivial $T$-Killing spinor whose energy-momentum tensor
$T^\psi$ corresponds to the orthogonal projection on
$\mathcal{K}^\perp$. Then the generalized cone $\mathcal{M}$ has to be
Ricci flat. 
\end{thm}
The argument used by C.~B\"ar in \cite{Bea} seems to be still usefull in
this situation. Therefore, we can formulate the following question in the
case where $(M^n,g)$ is a compact Riemannian spin manifold admitting a smooth parallel distribution.
\vspace{2mm}

\emph{Is it possible to find lower bounds for the first eigenvalue of the
  Dirac operator on $(M^n,g)$, such that the limiting cases are
  characterized by the existence of a non-trivial $T$-Killing spinor
  whose energy-momentum tensor is a projector?}

\vspace{2mm}

Such an estimate on a compact Riemannian spin manifold admitting a parallel
1-form is given in \cite{AGI}. One can note that the problem formulated
above is a particular case of the problem of studying the spectrum of the
Dirac operator on a manifold admitting a parallel $k$-form.

\providecommand{\bysame}{\leavevmode\hbox to3em{\hrulefill}\thinspace}

\end{document}